\documentclass[11pt]{amsart}
\usepackage[T1]{fontenc}
\usepackage{lmodern,amssymb,mathtools,microtype,booktabs,needspace,graphicx,array}
\usepackage[letterpaper,margin=1.08in]{geometry}
\usepackage{xcolor}
\usepackage[colorlinks=true,linkcolor=blue,citecolor=blue,urlcolor=blue]{hyperref}
\hypersetup{pdftitle={On dual varieties of smooth projective varieties},pdfauthor={Timofey Fedorov},pdfsubject={Projective duality}}
\newtheorem{theorem}{Theorem}[section]
\newtheorem{proposition}[theorem]{Proposition}
\newtheorem{lemma}[theorem]{Lemma}
\newtheorem{corollary}[theorem]{Corollary}
\newtheorem{maintheorem}{Theorem}
\newtheorem{maincorollary}[maintheorem]{Corollary}

\newcommand{\CC}{\mathbb{C}}
\newcommand{\PP}{\mathbb{P}}
\newcommand{\OO}{\mathcal O}
\numberwithin{equation}{section}

\title[On dual varieties of smooth projective varieties]{On dual varieties of smooth projective varieties}
\author{Timofey Fedorov}
\address{National Research University Higher School of Economics, Moscow, \mbox{Russia}.}
\email{tsfedorov@edu.hse.ru}
\date{}
\begin{document}
\raggedbottom
\begin{abstract}
For a smooth irreducible nondegenerate complex projective
variety $X\subsetneq\PP^N$, we prove that its projective dual
$X^\vee$ is a hypersurface whenever
$\dim X\ge3\operatorname{codim}X-2$.
This result may be regarded as further evidence for
Hartshorne's complete intersection conjecture.
We also prove the duality defect conjecture of Landman
in codimension at most five. Finally, we complete Ein's classification of smooth nondegenerate
complex projective varieties with smooth projective dual.
The list consists of smooth quadrics, the Segre varieties
$\PP^1\times\PP^m$ with $m\ge2$, the Grassmannian
$\operatorname{Gr}(2,5)$, and the spinor tenfold in $\PP^{15}$.
\end{abstract}

\maketitle

\section*{Introduction}

We work over $\CC$. Let $X\subsetneq\PP^N$ be a smooth
nondegenerate projective variety, and let
$X^\vee$ be its projective dual. Write
$k=\operatorname{codim} X^\vee-1$ for the duality defect of $X$. The duality defect conjecture predicts that $X^\vee$ is a
hypersurface whenever $\dim X>2\operatorname{codim}X$.
It would follow from Hartshorne's complete intersection
conjecture: in this range, $X$ would be a complete
intersection, whose projective dual is a hypersurface
\cite[Proposition~3.1]{Ein}.

The duality defect conjecture was proved in codimension
two by Ein \cite[Theorem~3.4]{Ein}.
Jorgenson \cite{Jorgenson} proved it in codimension three
when the dimension of the projective space is odd and
established it in certain finite ranges in codimensions
three, four, and five. In particular, his results cover
$\operatorname{codim}X=4$ for $10\le\dim X\le46$
and $\operatorname{codim}X=5$ for $13\le\dim X\le18$.

For varieties with positive duality defect, we prove the following
lower bound for duality defect, announced by Zak \cite{Zak1991}.

\begin{maintheorem}
\label{dd:cor:zak-bound}
\label{dd:thm:range}
\label{dd:vanishing}
Let $X\subsetneq\PP^N$ be a smooth irreducible nondegenerate
complex projective variety, and let
$k=\operatorname{codim}X^\vee-1$ be its duality defect.
If $k>0$, then
\[
 k\ge\dim X-2\operatorname{codim}X+2.
\]
Consequently, $X^\vee$ is a hypersurface whenever
$\dim X\ge3\operatorname{codim}X-2$.
\end{maintheorem}

The proof uses two results of Biran and Jerby
\cite[Theorem~A and Corollary~B]{BiranJerby}:
the invertibility of the hyperplane class in the integral
quantum cohomology of a smooth hyperplane section,
and two-periodicity of certain sums of Betti numbers.
Combined with the integral Barth--Larsen theorem
\cite{Larsen}, these results allow us to determine
the Euler characteristics of $X$ and a general hyperplane
section and establish restrictions on the Chern classes
of the normal bundle of $X$. These Chern class restrictions imply that both Euler
characteristics are divisible by $\deg X$, yielding
the lower bound on the duality defect.

For $1\le\operatorname{codim}X\le5$ and
$\dim X>2\operatorname{codim}X$, the only pairs not covered
by our Theorem~\ref{dd:cor:zak-bound} are
\[
 (9,4),\qquad (11,5),\qquad (12,5).
\]
The first two pairs are excluded in
Propositions~\ref{dd:prop:nine-four} and
\ref{dd:prop:eleven-five}, respectively.
The exclusion of $(12,5)$ follows from that of $(11,5)$
by taking a general hyperplane section.

\begin{maincorollary}\label{dd:cor:codimension-five}
Let $X\subsetneq\PP^N$ be a smooth irreducible nondegenerate
complex projective variety with
$1\le\operatorname{codim}X\le5$.
If
\[
 \dim X>2\operatorname{codim}X,
\]
then $X^\vee$ is a hypersurface.
\end{maincorollary}

For completeness, Proposition~\ref{dd:small-codim}
lists the pairs $(\dim X,\operatorname{codim}X)$ with
$1\le\operatorname{codim}X\le5$ that admit smooth
nondegenerate projective varieties with positive duality defect.
Its proof uses Corollary~\ref{dd:cor:codimension-five}
and classification results of Ein, Ionescu and Russo,
and Bayer--Kuznetsov--Macr\`i \cite{Ein,IRdef,BKM}.
Each pair in the list is realized by a Segre variety,
the Grassmannian $\operatorname{Gr}(2,5)$,
the spinor tenfold, or a suitable linear section
of one of these varieties.

Finally, we complete the classification of smooth
nondegenerate projective varieties with smooth projective dual.
Ein obtained the classification under the assumption
$\dim X\le2\operatorname{codim} X$ \cite[Theorem~4.5]{Ein}.
We prove that, when $\dim X>2\operatorname{codim} X$, the only
such varieties are quadric hypersurfaces. 

\begin{maintheorem}\label{hd:dual:classification}
Let $X\subsetneq\PP^N$ be a smooth irreducible nondegenerate
complex projective variety. Then
$X^\vee$ is smooth if and only if $X$
is projectively equivalent to one of the following varieties:
\begin{enumerate}
\item a smooth quadric hypersurface
      $Q^{\dim X}\subset\PP^{\dim X+1}$;
\item the Segre variety
      $\PP^1\times\PP^{\dim X-1}\subset\PP^{2\dim X-1}$,
      with $\dim X\ge3$;
\item the Pl\"ucker Grassmannian
      $\operatorname{Gr}(2,5)\subset\PP^9$;
\item the spinor tenfold
      $\operatorname{OG}^+(5,10)\subset\PP^{15}$
      in its half-spin embedding.
\end{enumerate}
\end{maintheorem}
The main step in the proof of
Theorem~\ref{hd:dual:classification} is
Proposition~\ref{hd:dual:exclusion}, which shows that
$X^\vee$ is singular whenever
$\operatorname{codim}X\ge4$ and
$\dim X-\operatorname{codim}X\ge6$.
Its proof combines a Chern character calculation on
the conormal variety (Lemma~\ref{hd:dual:chern})
with Hwang's nondegeneracy theorem
\cite[Corollary~2.8]{HwangLines}.

The paper is organized as follows. In
Section~\ref{sec:duality}, we prove
Theorem~\ref{dd:cor:zak-bound} and
Corollary~\ref{dd:cor:codimension-five}.
In Section~\ref{hd:dual:smooth}, we prove
Theorem~\ref{hd:dual:classification}.

\subsection*{Notation and conventions}

Let $X\subset\PP^N$ be a smooth irreducible nondegenerate
complex projective variety of dimension $n$ and codimension
$c\ge1$, and let $Y$ be a general hyperplane section of $X$.
We write $d=\deg X$, $h=c_1(\OO_X(1))$, and
$\mathcal N=N_{X/\PP^N}$.

The projective dual $X^\vee\subset(\PP^N)^\vee$ is the closure
of the set of hyperplanes tangent to $X$.
For a smooth nondegenerate projective variety $V$, write
$k(V)=\operatorname{codim}V^\vee-1$ for its duality defect.
The variety $V$ is called dual defective if $k(V)>0$.
For $X$ as above, put $k=k(X)$.

For a smooth projective variety $V$, we write
\[
 H^*(V)=H^*(V,\mathbb Z)/\text{torsion}.
\]
We use the same notation for an integral cohomology class
and its image in this quotient, and write $b_j(V)$
for the $j$th Betti number of $V$.

We write $\chi(V)$ for the topological Euler characteristic
of $V$, and $\chi(V,F)$ for the Euler characteristic of a
coherent sheaf $F$ on $V$.

For a vector bundle $F$ on $V$, we regard its Chern classes
$c_j(F)$ as elements of $H^{2j}(V)$, retaining integral
coefficients for the divisibility arguments in
Section~\ref{sec:duality}. In the Chern character calculations
of Section~\ref{hd:dual:smooth}, we work in rational cohomology,
regarding the Chern classes as elements of
$H^{2j}(V,\mathbb Q)$, where the components
$\operatorname{ch}_j(F)$ of the Chern character also lie.

We use the convention that $\PP(F)$ parametrizes
one-dimensional linear subspaces of the fibres of $F$,
with tautological line bundle $\OO_{\PP(F)}(-1)$.

For quantum cohomology we follow \cite{BiranJerby}.
In the situations considered below, $QH^*(Y;\Lambda)$
denotes the small quantum cohomology ring with
$\Lambda=\mathbb Z[q,q^{-1}]$ and $\deg q=n+k$.
Its underlying graded abelian group is
$H^*(Y)\otimes_{\mathbb Z}\Lambda$.
We denote multiplication in this ring by $*$,
to distinguish it from ordinary cup product.

\medskip
\noindent\textbf{AI usage.}
The author used ChatGPT (OpenAI) to help develop and draft mathematical
arguments and examples, search the literature, edit the prose, and
prepare the \LaTeX{} source. The author is responsible for the final text,
including its mathematical claims, proofs, and references.

\medskip
\noindent\textbf{Acknowledgements.}
I am grateful to my advisor Serge Lvovski for his support and guidance.
The study has been funded within the framework of the HSE University Basic
Research Program (HSE-BR-2025-060).

\Needspace{10\baselineskip}
\section{Duality defect}\label{sec:duality}

\subsection{Euler characteristic of dual defective varieties}
\label{dd:subsec:small-codimension}

\begin{proposition}\label{dd:prop:euler}
Let $X\subsetneq\PP^N$ be a smooth irreducible nondegenerate
complex projective variety, and let
$Y$ be a general hyperplane section of $X$.
Put $k=\operatorname{codim}X^\vee-1$. Suppose that
$\dim X>2\operatorname{codim}X$ and $k>0$.
Set $\ell_X=\dim X+k+2$ and $\ell_Y=\dim X+k$.
For $V=X,Y$ and $j\in\mathbb Z$, set
\[
 \widetilde b_j(V)=\sum_{t\in\mathbb Z} b_{j+t\ell_V}(V).
\]
Then
\begin{equation}\label{dd:eq:folded-betti}
 \widetilde b_j(X)=\widetilde b_j(Y)=
 \begin{cases}2,&j\text{ even},\\0,&j\text{ odd}.\end{cases}
\end{equation}
In particular,
\begin{equation}\label{dd:eq:euler}
 \chi(X)=\dim X+k+2,\qquad \chi(Y)=\dim X+k.
\end{equation}
\end{proposition}

\begin{proof}
Zak's upper bound on the duality defect
\cite[Ch.~I, Corollary~2.5]{Zak} gives $k\le c-1$.
Landman's parity theorem
\cite[Theorem~2.4(a)]{Ein} gives $k\equiv n\pmod2$. Write $Y=X\cap L$, where $L\subset\PP^N$ is a general hyperplane.
Both $X\subset\PP^N$ and $Y\subset L\simeq\PP^{N-1}$ have
codimension $c$. By the Barth--Larsen theorem \cite{Larsen}, the
restriction maps
\begin{align*}
 H^j(\PP^N,\mathbb Z)&\longrightarrow H^j(X,\mathbb Z)
 &&(0\le j\le n-c),\\
 H^j(L,\mathbb Z)&\longrightarrow H^j(Y,\mathbb Z)
 &&(0\le j\le n-c-1)
\end{align*}
are isomorphisms in the indicated ranges. In particular, $b_2(X)=1$. Since $k>0$, Corollary~B of
\cite{BiranJerby} gives $\widetilde b_{j+2}(V)=\widetilde b_j(V)$
for both $X$ and its smooth hyperplane section $Y$.

Thus, it suffices to compute $\widetilde b_j(V)$
in two consecutive degrees. Put
$r=(n-k-2)/2$. The inequalities $n>2c$ and
$1\le k\le c-1$ give $1\le r$ and $r+1\le n-c-1$.
For $V=X,Y$, we have $\ell_V=2\dim V-2r$.
Thus, for $j=r,r+1$, only the summands indexed by
$t=0,1$ in the definition of $\widetilde b_j(V)$
can be nonzero. By Poincar\'e duality,
\[
 \widetilde b_r(V)=2b_r(V),\qquad
 \widetilde b_{r+1}(V)=b_{r+1}(V)+b_{r-1}(V).
\]
The three degrees $r-1,r,r+1$ all lie in the Barth--Larsen range
for both varieties. Their Betti numbers are therefore one in even
degrees and zero in odd degrees. Two-periodicity gives
\eqref{dd:eq:folded-betti}. By Landman's parity theorem,
$n+k$ is even, so $\ell_X=n+k+2$ and $\ell_Y=n+k$
are both even. Using \eqref{dd:eq:folded-betti}, we obtain,
for $V=X,Y$,
\[
 \chi(V)=\sum_{j=0}^{\ell_V-1}
 (-1)^j\widetilde b_j(V)=\ell_V.
\]
\end{proof}

\begin{lemma}\label{dd:lem:divisibility}
Let $X\subsetneq\PP^N$ be a smooth irreducible nondegenerate
projective variety, and let $Y$ be a
general hyperplane section of $X$. Put $h=c_1(\OO_X(1))$
and $\mathcal N=N_{X/\PP^N}$. Suppose that
\begin{equation}\label{dd:eq:chern-integral}
 c_j(\mathcal N)\in\mathbb Z h^j
 \quad\text{in }H^{2j}(X),\qquad
 1\le j<\operatorname{codim}X.
\end{equation}
Then $\deg X$ divides both $\chi(X)$ and $\chi(Y)$.
Condition~\eqref{dd:eq:chern-integral} holds whenever
$\dim X\ge3\operatorname{codim}X-2$.
\end{lemma}

\begin{proof}
Let $i:X\hookrightarrow\PP^N$ be the inclusion. The self-intersection
formula gives $c_c(\mathcal N)=i^*i_*(1)=d h^c$. Together with
\eqref{dd:eq:chern-integral}, this shows that $c(\mathcal N)$ lies in
the subring $\mathbb Z[h]\subset H^*(X)$. Its inverse
also lies in this subring, since $c(\mathcal N)-1$ is nilpotent.

The Euler sequence on $\PP^N$ and the normal sequence of $X$ give
\[
 c(T_X)=(1+h)^{N+1}c(\mathcal N)^{-1}\in\mathbb Z[h].
\]
Put $h_Y=h|_Y$. Since $N_{Y/X}\simeq\OO_Y(1)$, the normal sequence
of $Y\subset X$ gives
\[
 c(T_Y)=c(T_X)|_Y\,(1+h_Y)^{-1}\in\mathbb Z[h_Y].
\]
Thus, modulo torsion, $c_n(T_X)=u h^n$ and
$c_{n-1}(T_Y)=v h_Y^{n-1}$ for some integers $u,v$.
We obtain
\[
 \chi(X)=\int_X c_n(T_X)=ud,\qquad
 \chi(Y)=\int_Y c_{n-1}(T_Y)=vd.
\]
Hence $d$ divides both $\chi(X)$ and $\chi(Y)$.

If $n\ge3c-2$, then $2j\le2(c-1)\le n-c$ for $1\le j<c$.
Barth--Larsen therefore gives $H^{2j}(X,\mathbb Z)=\mathbb Z h^j$
in this range, proving \eqref{dd:eq:chern-integral}.
\end{proof}

\subsection{A lower bound for positive duality defect}\label{dd:subsec:integral}

\begin{proposition}\label{dd:prop:primitive-powers}
Let $X\subsetneq\PP^N$ be a smooth irreducible nondegenerate
complex projective variety. Put
$h=c_1(\OO_X(1))$ and $k=\operatorname{codim}X^\vee-1$.
If $\dim X>2\operatorname{codim}X$ and $k>0$, then
\begin{equation}\label{dd:eq:primitive-powers}
 H^{2j}(X)=\mathbb Z h^j
 \qquad\text{for }0\le2j\le\dim X-k-2.
\end{equation}
In particular,
\begin{equation}\label{dd:eq:defect-lower-bound}
 k\ge\dim X-2\operatorname{codim}X+2.
\end{equation}
\end{proposition}

\begin{proof}
Put $h_Y=h|_Y$, $m=n-1$, $C=(n+k)/2$, and $r=m-C=(n-k-2)/2$.
By Landman's parity theorem, $C$ and $r$ are integers. Since $n>2c$
and $1\le k\le c-1$, we have $0\le r<C$. Since $k>0$ and $b_2(X)=1$,
\cite[Theorem~A]{BiranJerby} implies that $h_Y$
is invertible in $QH^*(Y;\Lambda)$. The underlying
graded abelian group of $QH^*(Y;\Lambda)$ is
$H^*(Y)\otimes_{\mathbb Z}\Lambda$, where
$\Lambda=\mathbb Z[q,q^{-1}]$ and $\deg q=2C$.

For $0\le j\le r$, we have $j<C$ and $j+C\le m<2C$. Hence
\begin{equation}\label{dd:eq:two-summands}
 QH^{2j}(Y;\Lambda)
 =H^{2j}(Y)\oplus H^{2j+2C}(Y)q^{-1}.
\end{equation}
By \eqref{dd:eq:folded-betti}, this group has rank two.
Since $0\le j\le j+C\le m$ and
$\int_Y h_Y^m=d>0$, both $h_Y^j$ and $h_Y^{j+C}$
are nonzero. Hence each summand in
\eqref{dd:eq:two-summands} has rank one. In particular,
$H^{2j}(Y)\simeq\mathbb Z$.

It remains to show that $h_Y^j$ generates $H^{2j}(Y)$. Write
$h_Y^{*j}$ for the $j$th power of $h_Y$ in
$QH^*(Y;\Lambda)$. For $0\le j<C$, we have
$h_Y^{*j}=h_Y^j$. Indeed, the term independent of $q$
is the ordinary cup-product power, while any term involving
$q^\ell$ with $\ell>0$ has coefficient in
$H^{2j-2C\ell}(Y)$. Since $2j-2C\ell<0$,
no such terms can occur. The class $1$ generates the direct summand
$H^0(Y)=\mathbb Z$ of $QH^0(Y;\Lambda)$.
Since $h_Y$ is invertible, multiplication by
$h_Y^{*j}$ induces an isomorphism
\[
 QH^0(Y;\Lambda)\xrightarrow{\ \sim\ }QH^{2j}(Y;\Lambda),
 \qquad 1\longmapsto h_Y^{*j}.
\]
Consequently, $\mathbb Z h_Y^j$ is a direct summand
of $QH^{2j}(Y;\Lambda)$. For $0\le j\le r$,
the class $h_Y^j$ belongs to the rank-one summand
$H^{2j}(Y)$ in \eqref{dd:eq:two-summands}
and therefore generates this group.
Since $2j\le2r<n-1=\dim Y$, the integral weak
Lefschetz theorem gives an isomorphism
\[
 H^{2j}(X,\mathbb Z)\xrightarrow{\ \sim\ }
 H^{2j}(Y,\mathbb Z)
\]
induced by restriction, which sends $h^j$ to $h_Y^j$.
Passing to cohomology modulo torsion, we obtain
$H^{2j}(X)=\mathbb Z h^j$, proving
\eqref{dd:eq:primitive-powers}.

Suppose that $k\le n-2c$. Then
$2(c-1)\le n-k-2$, so
\eqref{dd:eq:primitive-powers} implies that
\[
 c_j(\mathcal N)\in\mathbb Z h^j
 \qquad\text{for }1\le j<c.
\]
By Lemma~\ref{dd:lem:divisibility} and
Proposition~\ref{dd:prop:euler}, we obtain
$d\mid\chi(X)-\chi(Y)=2$.
This contradicts $d\ge c+1\ge3$.
Hence $k>n-2c$, and Landman's parity theorem gives
\eqref{dd:eq:defect-lower-bound}.
\end{proof}

\begin{proof}[Proof of Theorem~\ref{dd:cor:zak-bound}]
Assume that $k>0$. If $n>2c$, the lower bound follows
from Proposition~\ref{dd:prop:primitive-powers}.
If $n<2c$, it follows from
$n-2c+2\le1\le k$.
Finally, if $n=2c$, Landman's parity theorem
\cite[Theorem~2.4(a)]{Ein} gives $k\ge2$.
Combining the lower bound with Zak's upper bound
\cite[Ch.~I, Corollary~2.5]{Zak}, we obtain
\[
 \dim X-2\operatorname{codim}X+2
 \le k\le\operatorname{codim}X-1
\]
whenever $k>0$. Hence
$\dim X\le3\operatorname{codim}X-3$.
Thus $X^\vee$ is a hypersurface whenever
$\dim X\ge3\operatorname{codim}X-2$.
\end{proof}

\begin{proposition}\label{dd:prop:nine-four}
Let $X\subset\PP^{13}$ be a smooth irreducible nondegenerate
complex projective variety of dimension nine.
Then $X^\vee$ is a hypersurface.
\end{proposition}

\begin{proof}
Suppose that $k>0$. Theorem~\ref{dd:cor:zak-bound}
gives $k\ge3$, while Zak's upper bound
\cite[Ch.~I, Corollary~2.5]{Zak} gives $k\le3$.
Hence $k=3$.
A general hyperplane section $Y\subset X$ is smooth
and has dimension eight, codimension four, and duality
defect two \cite{EinII}.

Put $H=\OO_Y(1)$. By the integral Barth--Larsen theorem
\cite{Larsen}, $\operatorname{Pic}(Y)=\mathbb Z[H]$.
Ein's index formula \cite[Theorem~2.4(d)]{Ein} gives
$-K_Y=6H$.
Since $k(Y)=\operatorname{codim}Y-2$ and
$\dim Y=2\operatorname{codim}Y$,
\cite[Proposition~10(vii)]{IRdef} implies that $Y$
is a hyperplane section of the Grassmannian
$\operatorname{Gr}(2,5)$ or of the spinor tenfold.
Such sections have dimensions five and nine,
respectively, contradicting $\dim Y=8$.
\end{proof}

\begin{corollary}\label{dd:cor:chern-obstruction}
Let $X\subsetneq\PP^N$ be a smooth irreducible nondegenerate
complex projective variety. Put
$h=c_1(\OO_X(1))$, $\mathcal N=N_{X/\PP^N}$, and
$k=\operatorname{codim}X^\vee-1$.
If $\dim X>2\operatorname{codim}X$ and $k>0$,
then there is an integer $j$ such that
\[
 \frac{\dim X-\operatorname{codim}X}{2}
 <j<\operatorname{codim}X,
 \qquad c_j(\mathcal N)\notin\mathbb Z h^j
 \quad\text{in }H^{2j}(X).
\]
\end{corollary}

\begin{proof}
Suppose that
$c_j(\mathcal N)\in\mathbb Z h^j$ for every integer
$j$ satisfying $(n-c)/2<j<c$.
By Barth--Larsen, the same inclusion holds for
$1\le j<c$ with $2j\le n-c$.
Thus condition~\eqref{dd:eq:chern-integral} holds,
so Lemma~\ref{dd:lem:divisibility} and
Proposition~\ref{dd:prop:euler} give
$d\mid\chi(X)-\chi(Y)=2$, contradicting
$d\ge c+1\ge3$.
\end{proof}

\subsection{The case of dimension eleven and codimension five}
\label{dd:subsec:eleven-five}

\begin{proposition}\label{dd:prop:eleven-five}
Let $X\subset\PP^{16}$ be a smooth irreducible nondegenerate
complex projective variety with $\dim X=11$.
Then $X^\vee$ is a hypersurface.
\end{proposition}

\begin{proof}
Suppose, for a contradiction, that the duality defect
$k=\operatorname{codim}X^\vee-1$ is positive.
Since $\operatorname{codim}X=5$, Theorem~\ref{dd:cor:zak-bound}
gives $k\ge3$. Zak's theorem on tangencies gives $k\le4$
\cite[Ch.~I, Corollary~2.5]{Zak}, while Landman's parity theorem
\cite[Theorem~2.4(a)]{Ein} gives $k\equiv11\pmod2$.
Hence $k=3$.

Put $H=\OO_X(1)$ and $h=c_1(H)$. Since
$\dim X-\operatorname{codim}X=6$, the integral
Barth--Larsen theorem \cite{Larsen} gives
\[
 H^j(X,\mathbb Z)=
 \begin{cases}
  \mathbb Z h^{j/2},&j\text{ even},\\
  0,&j\text{ odd},
 \end{cases}
 \qquad 0\le j\le6.
\]
Hodge decomposition and the exponential sequence then give
$\operatorname{Pic}(X)=\mathbb Z[H]$.
Ein's index formula \cite[Theorem~2.4(d)]{Ein} yields
\[
 \OO_X(-K_X)\simeq H^{\otimes8}.
\]
Thus $X$ is Fano of index eight.

Since $3\dim X=33>2(16-1)=30$, Zak's linear normality
theorem \cite[Ch.~II, Corollary~2.17]{Zak} shows that
$X\subset\PP^{16}$ is linearly normal. As $X$ is also
nondegenerate, the restriction map
\[
 H^0(\PP^{16},\OO_{\PP^{16}}(1))
 \longrightarrow H^0(X,H)
\]
is an isomorphism. Hence $h^0(X,H)=17$.

Let $d=\deg X$ and $P(t)=\chi(X,H^{\otimes t})$.
Kodaira vanishing gives $H^i(X,H^{\otimes t})=0$ for
$i>0$ and every integer $t\ge-7$. Since negative powers
of $H$ have no nonzero global sections, we obtain
\[
 P(-1)=\cdots=P(-7)=0,\qquad P(0)=1,\qquad P(1)=17.
\]
Serre duality gives $P(-8-t)=-P(t)$. Moreover, $P$ has
degree $11$ and leading coefficient $d/11!$.

The quotient $11!P(t)/\prod_{j=1}^{7}(t+j)$ is a polynomial
of degree four with leading coefficient $d$. Since it is
invariant under $t\mapsto-8-t$, it is an even polynomial
in $t+4$. The values $P(0)=1$ and $P(1)=17$ determine
the remaining two coefficients, giving
\[
 P(t)=\frac{\prod_{j=1}^{7}(t+j)}{11!}
 \bigl(d(t+4)^4+(990-41d)(t+4)^2+400d-7920\bigr).
\]

Since $H^4(X,\mathbb Z)=\mathbb Z h^2$, write
$c_2(T_X)=b h^2$ for some $b\in\mathbb Z$. As
$c_1(T_X)=8h$, comparison of the coefficient of $t^9$
in the preceding formula with Hirzebruch--Riemann--Roch gives
\[
 \frac{825d+990}{11!}
 =\frac{(64+b)d}{12\cdot9!}.
\]
Hence $(b-26)d=108$, and therefore $d\mid108$.

Since $X\subset\PP^{16}$ is nondegenerate, the degree bound
$\deg X\ge\operatorname{codim}X+1$ gives $d\ge6$.
If $3\nmid d$, then $d\mid108=2^2\cdot3^3$ implies
$d\mid4$, a contradiction. Thus $3\mid d$.

Let $Y\subset X$ be a smooth hyperplane section and put
$a=h|_Y$. Then $\dim Y=10$ and $\int_Y a^{10}=d$.
Integral weak Lefschetz identifies $H^6(X,\mathbb Z)$
with $H^6(Y,\mathbb Z)$ by restriction. Hence
\[
 H^6(Y,\mathbb Z)=\mathbb Z a^3.
\]

Since $k=3>0$ and $b_2(X)=1$,
\cite[Theorem~A]{BiranJerby} shows that $a$ is
invertible in $QH^*(Y;\Lambda)$. The underlying graded abelian group of
$QH^*(Y;\Lambda)$ is
$H^*(Y)\otimes_{\mathbb Z}\Lambda$, where
\[
 \Lambda=\mathbb Z[q,q^{-1}],\qquad
 \deg q=\dim X+k=14.
\]

Since $\dim Y=10$ and $\deg q=14$, the graded components
in degrees $6$, $8$, and $10$ are
\[
\begin{aligned}
 QH^6(Y;\Lambda)
   &=H^6(Y)\oplus H^{20}(Y)q^{-1}
     \simeq\mathbb Z^2,\\
 QH^8(Y;\Lambda)&=H^8(Y),\\
 QH^{10}(Y;\Lambda)&=H^{10}(Y).
\end{aligned}
\]
Since $a$ is invertible in $QH^*(Y;\Lambda)$,
multiplication by $a$ in this ring induces isomorphisms
\[
 QH^6(Y;\Lambda)
 \xrightarrow{\ a*\ } QH^8(Y;\Lambda)
 \xrightarrow{\ a*\ } QH^{10}(Y;\Lambda).
\]
Consequently,
\[
 \operatorname{rank}H^8(Y)
 =\operatorname{rank}H^{10}(Y)=2.
\]

Since the products $a*a^3$ and $a*a^4$ have degrees
$8$ and $10$, both smaller than $\deg q=14$,
no terms involving positive powers of $q$ can occur.
Thus
\[
 a*a^3=a^4,\qquad a*a^4=a^5.
\]
In the decomposition of $QH^6(Y;\Lambda)$ above,
$a^3$ generates the direct summand
$H^6(Y)\simeq\mathbb Z$.
The composition of the two isomorphisms above
sends $a^3$ to $a^5$. Consequently,
$\mathbb Z a^5$ is a direct summand of $H^{10}(Y)$.

Since $\dim_{\mathbb C}H^8(Y,\mathbb C)=2$ and
$a^4\ne0$ has Hodge type $(4,4)$, Hodge symmetry
rules out any component of type $(p,q)$ with $p\ne q$.
Thus $H^8(Y,\mathbb C)=H^{4,4}(Y)$.
By the Lefschetz decomposition,
\[
 H^8(Y,\mathbb R)
 =\mathbb R a^4\oplus P^8(Y,\mathbb R),
\]
where $P^8(Y,\mathbb R)$ is the one-dimensional space
of primitive cohomology classes, defined by
\[
 P^8(Y,\mathbb R)=
 \ker\bigl(a^3\smile-:
 H^8(Y,\mathbb R)\longrightarrow H^{14}(Y,\mathbb R)\bigr).
\]
Since $\deg q=14$,
multiplication by $a$ in $QH^*(Y;\Lambda)$, from degree
$8$ to degree $10$, agrees with the ordinary cup-product map
\[
 a\smile-:H^8(Y)\longrightarrow
 H^{10}(Y).
\]
This map is therefore an isomorphism, and
\[
 H^{10}(Y,\mathbb R)
 =\mathbb R a^5\oplus aP^8(Y,\mathbb R).
\]

The two summands in this decomposition are orthogonal
with respect to the intersection pairing. Indeed, if
$\alpha\in P^8(Y,\mathbb R)$, then $a^3\alpha=0$ by
definition, so
\[
 \int_Y a^5\smile(a\alpha)=\int_Y a^6\alpha=0.
\]
Moreover,
\[
 \int_Y(a^5)^2=d>0.
\]
Since every class in $P^8(Y,\mathbb R)$ has Hodge type
$(4,4)$, the Hodge--Riemann bilinear relations give
\[
 \int_Y a^2\alpha^2>0
 \qquad\text{for every }0\ne\alpha\in P^8(Y,\mathbb R).
\]
Thus the intersection pairing on $H^{10}(Y,\mathbb R)$
is positive definite.

By Poincar\'e duality, the intersection pairing on
$H^{10}(Y)$ is integral and unimodular.
Since $\mathbb Z a^5$ is a direct summand of the
rank-two free abelian group $H^{10}(Y)$, we can
choose $w\in H^{10}(Y)$ such that $(a^5,w)$
is an integral basis. In this basis,
the intersection matrix has the form
\[
 \begin{pmatrix}
 d&e\\
 e&f
 \end{pmatrix},
 \qquad e,f\in\mathbb Z.
\]
Since the pairing is positive definite and unimodular,
its determinant is one. Hence
\[
 df-e^2=1.
\]
But $3\mid d$, so reduction modulo $3$ gives
$e^2\equiv-1\pmod3$, a contradiction. Therefore
$X^\vee$ is a hypersurface.
\end{proof}

\begin{proof}[Proof of Corollary~\ref{dd:cor:codimension-five}]
Put $n=\dim X$ and $c=\operatorname{codim}X$, and assume
$n>2c$. For $c\le3$, this implies $n\ge3c-2$, so
Theorem~\ref{dd:cor:zak-bound} applies. The same theorem
covers $c=4$ when $n\ge10$ and $c=5$ when $n\ge13$.
The remaining pairs are
\[
 (n,c)=(9,4),(11,5),(12,5).
\]
Propositions~\ref{dd:prop:nine-four} and
\ref{dd:prop:eleven-five} exclude the first two pairs.
For $(n,c)=(12,5)$, suppose that $k>0$.
Landman's parity theorem gives $k\ge2$.
A general hyperplane section $Y$ of $X$ is smooth
and has $\dim Y=11$, $\operatorname{codim}Y=5$,
and positive duality defect $k-1$ \cite{EinII}.
This contradicts Proposition~\ref{dd:prop:eleven-five}.
\end{proof}

\subsection{Dimensions of dual defective varieties in codimension at most five}
\label{hd:dual:smallcodim}

\begin{proposition}\label{hd:dual:smallcodim-table}\label{dd:small-codim}
The following table lists exactly those pairs
$(\dim X,\operatorname{codim}X)$ with
$1\le\operatorname{codim}X\le5$ for which there exists
a smooth irreducible nondegenerate complex projective variety
$X\subsetneq\PP^N$ with positive duality defect.
\begin{center}
\renewcommand{\arraystretch}{1.15}
\begin{tabular}{@{}cl@{}}
\toprule
$\operatorname{codim}X$ & Possible values of $\dim X$\\
\midrule
$1$ & none\\
$2$ & $3$\\
$3$ & $3,4,5,6$\\
$4$ & $3,4,5$\\
$5$ & $3,4,5,6,7,8,9,10$\\
\bottomrule
\end{tabular}
\end{center}
\end{proposition}

\begin{proof}
Let $X\subset\PP^{n+c}$ be a smooth irreducible nondegenerate
projective variety of dimension $n$ and codimension
$1\le c\le5$. If $n>2c$, Corollary~\ref{dd:cor:codimension-five}
gives $k(X)=0$. Hence it suffices to consider $1\le n\le2c$.

Curves, surfaces, and smooth complete intersections have defect zero
by \cite[Proposition~3.1]{Ein}; in particular, this excludes
$n\le2$ and $c=1$. For a variety with positive duality defect,
Zak's theorem on tangencies and Landman's parity theorem give
\begin{equation}\label{hd:dual:smallcodim-bounds}
 1\le k(X)\le c-1,\qquad k(X)\equiv n\pmod2;
\end{equation}
see \cite[Theorem~1.3 and Theorem~2.4(a)]{Ein}.
Consequently, $(n,c)=(4,2)$ is also excluded.
Among the pairs with $1\le n\le2c$ not listed in the table,
the only cases not yet excluded are
$(n,c)=(6,4),(7,4),(8,4)$.

Suppose that $c=4$, $n\in\{6,7,8\}$, and $k(X)>0$.
Since $n-c\ge2$, the Barth--Larsen theorem gives
$\operatorname{Pic}(X)=\mathbb Z[H]$, where $H=\OO_X(1)$
\cite{Larsen}. Ein's index formula
\cite[Theorem~2.4(d)]{Ein} then gives
\begin{equation}\label{hd:dual:smallcodim-index}
 -K_X=\frac{n+k(X)+2}{2}H.
\end{equation}
If $n\in\{6,8\}$, \eqref{hd:dual:smallcodim-bounds}
gives $k(X)=2=c-2$. Since $n\le2c$,
\cite[Proposition~10(vii)]{IRdef} implies that $X$
is a hyperplane section of
$\operatorname{Gr}(2,5)\subset\PP^9$ or of the spinor tenfold
$\operatorname{OG}^+(5,10)\subset\PP^{15}$.
Such sections have dimensions five and nine, respectively,
contradicting $n\in\{6,8\}$.

If $n=7$, \eqref{hd:dual:smallcodim-bounds} gives
$k(X)\in\{1,3\}$. If $k(X)=3=c-1$, then
$n\le2c+2$, and \cite[Proposition~10(vi)]{IRdef}
implies that $X$ is projectively equivalent to
$\operatorname{Gr}(2,5)$ or the spinor tenfold.
These varieties have dimensions six and ten, respectively,
contradicting $n=7$. Therefore $k(X)=1$, and
\eqref{hd:dual:smallcodim-index} gives $-K_X=5H$.

Let $C$ be the smooth curve cut out on $X$ by six general
hyperplanes, and let $g$ be its genus. By adjunction,
$K_C=H|_C$, so $2g-2=H^7$.
Since $X$ is nondegenerate of codimension four,
$H^7\ge5$ and hence $g\ge4$.
By \cite[Lemma~6.1(a)]{BKM},
\[
 h^0(X,H)=\dim X+g-1=g+6.
\]
Since $3\cdot7>2(11-1)$, Zak's linear normality theorem
\cite[Ch.~II, Corollary~2.17]{Zak} gives
$h^0(X,H)=12$, hence $g=6$.
Now \cite[Theorem~1.1(b)]{BKM} gives
$\dim X\le6$, contradicting $\dim X=7$.

For the examples, we use the hyperplane-section formula
\cite{EinII}: if $V\subset\PP^s$ is smooth with positive
duality defect and $L\subset\PP^s$ is a general hyperplane, then
\[
 k(V\cap L)=k(V)-1.
\]
For $2\le c\le5$, consider the Segre variety
$S=\PP^1\times\PP^c\subset\PP^{2c+1}$.
Its duality defect is $k(S)=c-1$.

For $3\le n\le c+1$, let $S_n$ be a general linear
section of $S$ of codimension $c+1-n$. Then
$S_n\subset\PP^{n+c}$ is a smooth nondegenerate
$n$-fold, and the hyperplane-section formula gives
\[
 k(S_n)=(c-1)-(c+1-n)=n-2>0.
\]

The Pl\"ucker Grassmannian
$\operatorname{Gr}(2,5)\subset\PP^9$ and the spinor tenfold
$\operatorname{OG}^+(5,10)\subset\PP^{15}$ are projectively
self-dual \cite[\S I.4, Table~1]{BucArxiv}.
Since their dimensions are six and ten, respectively,
their duality defects are two and four. The remaining pairs
in the table are realized by general linear sections $X$
of these varieties. For sections of codimension $r$ in
the Grassmannian and the spinor tenfold, respectively, we have
\[
\begin{aligned}
 (\dim X,\operatorname{codim}X,k(X))
   &=(6-r,3,2-r) &&(0\le r\le1),\\
 (\dim X,\operatorname{codim}X,k(X))
   &=(10-r,5,4-r) &&(0\le r\le3).
\end{aligned}
\]
\end{proof}

\section{Smooth varieties with smooth projective dual}\label{hd:dual:smooth}
Let $X\subsetneq\PP(W)=\PP^N$ be a smooth irreducible
nondegenerate projective variety.
Put $n=\dim X$ and $c=N-n$. In the Chern character calculations of this section, we
regard Chern classes and their pullbacks as classes in
rational cohomology.

\subsection{The conormal variety and a Chern character identity}
For this computation, suppose that $X^\vee$ is smooth.
The projective dual $X^\vee$ is nondegenerate since $X$ is
smooth and nondegenerate. Applying Zak's dimension inequality
\cite[Ch.~I, Corollary~2.5]{Zak} to $X$ and $X^\vee$,
and using biduality, we obtain
\[
 \dim X\le\dim X^\vee
 \le\dim (X^\vee)^\vee=\dim X.
\]
Hence $\dim X^\vee=n$ and
\[
 k(X)=k(X^\vee)=N-1-n=c-1.
\]

Let $\mathcal C\subset X\times X^\vee$ be the conormal variety. Its points
are the pairs $(z,[\xi])$ for which the hyperplane defined by $\xi$
is tangent to $X$ at $z$. Denote the projections by $p:\mathcal C\to X$ and
$q:\mathcal C\to X^\vee$, and put
\[
 E=N_{X/\PP^N}(-1),\qquad E'=N_{X^\vee/(\PP^N)^\vee}(-1).
\]
For $z=[w]\in X$, the fibre $E_z$ is naturally identified
with $W/T_w\widehat X$, where $\widehat X\subset W$ is the
affine cone over $X$. Thus $\PP(E_z^*)$ parametrizes the
hyperplanes tangent to $X$ at $z$. By biduality, $\mathcal C$ is
also the conormal variety of $X^\vee$. Hence
\[
 \mathcal C=\PP_X(E^*)=\PP_{X^\vee}((E')^*).
\]
In the projective bundle description $\mathcal C=\PP_X(E^*)$,
the tautological line at $(z,[\xi])$ is the line
$\CC\xi\subset W^*$. In the description
$\mathcal C=\PP_{X^\vee}((E')^*)$, it is the line
$\CC w\subset W$, where $z=[w]$. Consequently,
\[
 \OO_{\PP_X(E^*)}(1)\simeq q^*\OO_{X^\vee}(1),\qquad
 \OO_{\PP_{X^\vee}((E')^*)}(1)\simeq p^*\OO_X(1).
\]
Write
\[
 H=c_1(\OO_X(1)),\quad H'=c_1(\OO_{X^\vee}(1)),\qquad
 A=p^*H,\quad B=q^*H'.
\]

\begin{lemma}\label{hd:dual:chern}
Let $X\subsetneq\PP^N$ be a smooth irreducible nondegenerate
complex projective variety. Suppose
that $X^\vee$ is smooth,
$\operatorname{codim}X\ge4$, and
$\dim X-\operatorname{codim}X\ge6$.
Put
\[
 H=c_1(\OO_X(1)),\qquad E=N_{X/\PP^N}(-1),
\]
and
\[
 i=\frac{\dim X+\operatorname{codim}X+1}{2},
 \qquad
 a=\frac{\dim X-\operatorname{codim}X+1}{2}.
\]
Then $\operatorname{Pic}(X)=\mathbb Z[\OO_X(1)]$ and
\begin{equation}\label{hd:dual:identity}
 -K_X=iH,\qquad c_1(E)=aH,\qquad
 \operatorname{ch}_2(T_X)=\frac{i}{2}H^2.
\end{equation}
The last identity holds in $H^4(X,\mathbb Q)$.
\end{lemma}

\begin{proof}
By Barth--Larsen \cite{Larsen}, the restriction maps
\[
 H^j(\PP(W),\mathbb Z)\longrightarrow H^j(X,\mathbb Z),
 \qquad
 H^j(\PP(W^*),\mathbb Z)\longrightarrow H^j(X^\vee,\mathbb Z)
\]
are isomorphisms for $0\le j\le n-c$.
In particular, Hodge decomposition
and the exponential sequence give
$\operatorname{Pic}(X)=\mathbb Z[\OO_X(1)]$ and
$\operatorname{Pic}(X^\vee)=\mathbb Z[\OO_{X^\vee}(1)]$.
Both varieties have duality defect $k=c-1>0$, so Ein's index formula
\cite[Theorem~2.4(d)]{Ein} gives $-K_X=iH$ and
$-K_{X^\vee}=iH'$, with $i=(n+k+2)/2$.
Since $E=N_{X/\PP^N}(-1)$, adjunction gives
\[
 c_1(E)=c_1(N_{X/\PP^N})-cH
 =(N+1)H+K_X-cH
 =(n+1-i)H=aH.
\]
The same calculation for $X^\vee$ gives
$c_1(E')=aH'$.

The relative Euler sequences for the two projective bundle
structures give
\[
 \operatorname{ch}(T_{\mathcal C/X})
 =e^B p^*\operatorname{ch}(E^*)-1,
 \qquad
 \operatorname{ch}(T_{\mathcal C/X^\vee})
 =e^A q^*\operatorname{ch}((E')^*)-1.
\]
Combining these identities with the relative tangent
sequences for $p$ and $q$, we obtain
\begin{equation}\label{hd:dual:euler}
 p^*\operatorname{ch}(T_X)+e^B p^*\operatorname{ch}(E^*)
 =q^*\operatorname{ch}(T_{X^\vee})+e^A q^*\operatorname{ch}((E')^*).
\end{equation}
Since $n-c\ge6$, Barth--Larsen gives
$H^{2j}(X,\mathbb Q)=\mathbb Q H^j$ for $0\le j\le3$.
As $\operatorname{rk}E=c\ge4$, the projective bundle formula
for $p$ gives
\[
 H^6(\mathcal C,\mathbb Q)
 =\mathbb Q A^3\oplus\mathbb Q A^2B
  \oplus\mathbb Q AB^2\oplus\mathbb Q B^3.
\]
Write $\operatorname{ch}_2(E)=bH^2$, where $b\in\mathbb Q$.
Since $\operatorname{ch}_2(E^*)=\operatorname{ch}_2(E)$,
the coefficient of $A^2B$ on the left-hand side of
\eqref{hd:dual:euler} is $b$. On the right-hand side,
this coefficient is $-a/2$, since
$c_1((E')^*)=-aH'$. Thus $b=-a/2$.

Since $N_{X/\PP^N}=E(1)$, we have
\[
 \operatorname{ch}_2(N_{X/\PP^N})
 =\operatorname{ch}_2(E)+Hc_1(E)+\frac c2H^2
 =\left(b+a+\frac c2\right)H^2.
\]
The normal sequence then gives
\[
 \operatorname{ch}_2(T_X)
 =\left(\frac{N+1}{2}-b-a-\frac c2\right)H^2
 =\frac{n+1-a}{2}H^2
 =\frac i2H^2.\qedhere
\]
\end{proof}

\subsection{Lines through a general point}
\begin{proposition}\label{hd:dual:exclusion}
Let $X\subsetneq\PP^N$ be a smooth irreducible nondegenerate
complex projective variety. If
$\operatorname{codim}X\ge4$ and
$\dim X-\operatorname{codim}X\ge6$, then $X^\vee$ is singular.
\end{proposition}

\begin{proof}
Suppose that $X^\vee$ is smooth. By Lemma~\ref{hd:dual:chern},
$\operatorname{Pic}(X)=\mathbb Z[\OO_X(1)]$ and
$-K_X=iH$, where $i=(n+c+1)/2$.
Since $i>(n+1)/2$, $X$ is covered by lines.
For a general point $x\in X$, let
$M\subset\PP(T_xX)$ be the variety of tangent directions
to lines on $X$ passing through $x$.
It is smooth of pure dimension $r=i-2$
\cite[Lemma~2 and the discussion following it]{IRdef}. Since $2r=n+c-3\ge n-1=\dim\PP(T_xX)$,
any two irreducible components of $M$ must intersect.
As $M$ is smooth, distinct irreducible components
cannot meet. Hence $M$ is irreducible. By Hwang's nondegeneracy theorem,
$\langle M\rangle=\PP(T_xX)$ \cite[Corollary~2.8]{HwangLines}.

Let $\pi:U\to M$ be the universal family of lines, let
$e:U\to X$ be the evaluation map, and put
$L=c_1(\OO_M(1))$. By
\cite[Proposition~1.1]{AC},
\[
 -K_M\equiv\pi_*e^*\operatorname{ch}_2(T_X)+\frac r2L.
\]
Since every line parametrized by $M$ passes through $x$,
the universal family is noncanonically isomorphic to
\[
 U\simeq\PP_M(\OO_M\oplus\OO_M(-1)),
\]
where $\OO_M(-1)$ is the tautological line bundle on
$M\subset\PP(T_xX)$. Writing
$\xi=c_1(\OO_U(1))=e^*H$, the projective bundle formula gives
\[
 \pi_*e^*(H^2)=\pi_*(\xi^2)=L.
\]
Substituting \eqref{hd:dual:identity}, we obtain
\[
 -K_M\equiv\frac{i+r}{2}L=(r+1)L.
\]

Let $D$ be the smooth irreducible curve cut out on $M$ by $r-1$
general hyperplanes. By adjunction,
\[
 2g(D)-2=(K_M+(r-1)L)L^{r-1}=-2L^r.
\]
Since $g(D)\ge0$ and $L^r$ is a positive integer, we obtain
$L^r=1$. As an irreducible projective variety of degree one,
$M$ is a linear subspace of $\PP(T_xX)$. But
$r=i-2<n-1$, contradicting
$\langle M\rangle=\PP(T_xX)$.

\end{proof}

\subsection{Completion of the classification}
\begin{proof}[Proof of Theorem~\ref{hd:dual:classification}]
Suppose that $X^\vee$ is smooth. Then $\dim X^\vee=n$ and
$k(X)=c-1$. If $c=1$, then $X$ and $X^\vee$ are smooth hypersurfaces.
Write $e=\deg X$ and $e'=\deg X^\vee$.
Applying the degree formula for the dual of a smooth
hypersurface to both varieties gives
\[
 e'=e(e-1)^n,\qquad e=e'(e'-1)^n,
\]
where the second equality uses biduality.
Since $e\ge2$, the first equality gives $e'\ge e\ge2$.
The second gives $e\ge e'$. Thus $e=e'=2$, and $X$ is a quadric.

Assume $c\ge2$ and suppose that $n>2c$. If $c\le5$,
Corollary~\ref{dd:cor:codimension-five} gives $k(X)=0$,
contrary to $k(X)=c-1>0$. If $c\ge6$, then
$n-c\ge c+1\ge7$, so
Proposition~\ref{hd:dual:exclusion} implies that $X^\vee$
is singular, also a contradiction. Therefore $n\le2c$. Since $\dim X^\vee=n$, Ein's classification
\cite[Theorem~4.5]{Ein} shows that $X$ is projectively
equivalent to one of the varieties listed in items
(2)--(4) of Theorem~\ref{hd:dual:classification}.

Conversely, each variety in the list is projectively equivalent to
its dual \cite[\S I.4, Table~1]{BucArxiv}. Hence its dual is smooth.
\end{proof}

\bibliographystyle{amsplain}
\bibliography{projective_duality_6}
\end{document}